\documentclass[12pt,a4paper]{article}
\usepackage[T1]{fontenc}
\usepackage[utf8]{inputenc}
\usepackage{amsmath,amssymb,amsthm}
\usepackage[a4paper,margin=23mm]{geometry}
\usepackage[hidelinks]{hyperref}
\allowdisplaybreaks[2]
\newtheorem{theorem}{Theorem}[section]
\newtheorem{lemma}[theorem]{Lemma}
\newtheorem{transferlemma}{Lemma}[section]

\newtheorem{proposition}[theorem]{Proposition}
\newtheorem{corollary}[theorem]{Corollary}
\theoremstyle{remark}
\newtheorem{remark}[theorem]{Remark}
\newtheorem{problem}[theorem]{Problem}
\newcommand{\K}{\mathbb K}
\newcommand{\N}{\mathbb N}
\newcommand{\Lop}{\mathcal L}
\newcommand{\co}{\operatorname{co}}
\newcommand{\spn}{\operatorname{span}}
\newcommand{\cl}{\overline}
\newcommand{\rad}{\operatorname{r}}
\newcommand{\cnull}{c_0(\N,\K)}
\title{Hyperinvariant subspaces of compactoid operators\\ on non-Archimedean locally convex spaces}
\author{Wies\l aw \'Sliwa\\[3pt]
\small Institute of Mathematics, Faculty of Exact and Technical Sciences\\
\small University of Rzesz\'ow, Rzesz\'ow, Poland\\
\small \href{mailto:wsliwa@ur.edu.pl}{\texttt{wsliwa@ur.edu.pl}}}
\date{27 November 2026}
\hypersetup{pdftitle={Hyperinvariant subspaces of compactoid operators on non-Archimedean locally convex spaces},pdfauthor={Wies\l aw \'Sliwa},pdfsubject={Mathematics preprint},pdfkeywords={compactoid operators, hyperinvariant subspaces, non-Archimedean spaces}}
\begin{document}
\maketitle
\begin{abstract}
Over a complete nontrivially valued non-Archimedean field $\K$, the invariant and hyperinvariant subspace problems for nonzero compactoid operators on $c_0(\N,\K)$ are respectively equivalent to the corresponding problems on all infinite-dimensional sequentially complete Hausdorff locally convex spaces. Both problems are also equivalent to their respective semicompact and compact variants on arbitrary infinite-dimensional Hausdorff locally convex spaces. The reduction is based on a Banach factorization and a transfer argument for the full commutant. Over locally compact fields these problems have positive answers. We extend the quasi-null-vector theorem to sequentially complete locally convex spaces and use classical Fredholm theory to treat positive norm-growth radius over arbitrary complete fields. For operators commuting with a compactoid operator of positive radius, we characterize hyperinvariance through irreducible polynomial identities and the separating properties of the dual. We also show that closed hyperinvariant subspaces of operators on $c_0$ extend to hyperinvariant subspaces for the full commutant over arbitrary complete isometric valued extensions. We also give hyperinvariance criteria for small-range operators on spaces not of countable type. The general injective, dense-range, zero-radius compactoid case over non-locally compact fields is left open.
\end{abstract}


\section{Introduction}
The invariant subspace problem for non-Archimedean Banach spaces differs markedly from its classical counterpart. Every infinite-dimensional Banach space of countable type over a complete nontrivially valued non-Archimedean field admits a continuous linear operator with no nontrivial closed invariant subspace \cite{Sliwa}. Analogous negative results for certain K\"othe spaces were obtained by Kasprzak \cite[Theorems 3.8 and 4.1]{Kasprzak}; see also the corrigendum \cite{KasprzakCorr}. Compactoid operators, however, admit positive results under additional assumptions on the scalar field or on the operator.

Babahmed and El Asri \cite[Theorem 9]{BE} considered compact operators possessing a quasi-null vector over a locally compact field. El Asri, Kubzdela and Babahmed \cite[Theorem 3.5]{EKB} established the existence of a finite-dimensional hyperinvariant subspace for a compact operator of positive norm-growth radius on a Banach space with a basis. Their Theorem 3.8 gives a hyperinvariant subspace for nonzero compact operators over locally compact fields (the nonzero qualification is necessary). Their Theorems 3.10 and 3.14 concern, respectively, an algebraic obstruction to hyperinvariance and a quasi-null-vector condition for operators commuting with a compact operator. We extend these results in Proposition~\ref{prop:alg} and Theorem~\ref{thm:quasinull}, respectively; Theorem~\ref{thm:dichotomy} also answers their concluding question in the case where a commuting compactoid operator has positive radius.

Our principal contribution is the equivalence in Theorem~\ref{thm:equiv}, which reduces the invariant and hyperinvariant subspace problems for the operator classes considered here to the corresponding problems on $\cnull$. The essential step is Lemma~\ref{lem:transfer}: the full commutant is handled by composition with the two factorization maps, not by extending each commuting operator to a Banach completion. The full-commutant composition argument is classical in the Hilbert-space theory of quasi-similarity; see Hoover \cite[Theorem 2.1, p. 110]{Hoover}, who derives it from an argument of Sz.-Nagy and Foias. Here it is combined with a non-Archimedean Banach factorization to obtain the six equivalent assertions in Theorem~\ref{thm:equiv}, also for ordinary invariance. Thus no metrizability or basis assumption on the original locally convex space is needed. In the complex setting, Kramar \cite[Theorem 7]{Kramar97} studied invariant subspaces of operator algebras under compactness and range assumptions, and Kramar \cite[Theorem 2.5]{Kramar00} established density of operators with a hyperinvariant subspace in the algebra of locally bounded operators on complete locally convex spaces, for the topology of uniform convergence on bounded sets. Our reduction instead uses a Banach model of an individual operator and transfers hyperinvariance for its full commutant.

A second contribution is Theorem~\ref{thm:extension}, which preserves hyperinvariance under complete isometric extensions of the scalar field. In its proof, the relevant countable-type space is generated by the matrix entries of an individual commuting operator, rather than by the whole extension field. Corollary~\ref{cor:separablefield} consequently reduces the general problem to scalar fields separable in the valuation topology. Theorem~\ref{thm:quasinull} extends the quasi-null-vector method to sequentially complete locally convex spaces; Proposition~\ref{prop:alg} and Theorem~\ref{thm:dichotomy} give the algebraic criterion and the positive-radius commutant dichotomy.

The spectral ingredients in Section~4 are classical: Serre's Fredholm theory \cite{Serre}, Gruson's treatment of Banach and locally convex spaces \cite[Sections 4 and 5]{Gruson}, and the algebraic index-zero result for sequentially complete locally convex spaces \cite[Corollary 3, p. 116]{SliwaFred}. We use these ingredients to obtain the required locally convex formulation; the countable-type reduction of the positive-radius argument itself is already in \cite[Section 4.3, pp. 88--89]{Gruson}. The quasi-null-vector proof adapts \cite{BE,EKB}. Our proofs do not invoke the invariant or hyperinvariant subspace theorems of \cite{EKB} as black boxes.

In particular, the compactoid problem is settled over locally compact fields, in the positive norm-growth-radius case over arbitrary complete fields, and in the matrix-defined cases of Corollary~\ref{cor:extension}. The principal unresolved case concerns injective compactoid endomorphisms of $\cnull$ with dense range and zero norm-growth radius over non-locally compact fields. Section~\ref{sec:questions} states the remaining questions.

The compactoid case on spaces not of countable type, by contrast, follows directly from the classical countable-type theorem for the range of a compactoid operator \cite[Theorem 8.1.3(iv)]{PGS}; in the Banach case this observation was already made in \cite[Remark 10, pp. 212--213]{BE}. In the final section (Section~\ref{sec:cardinal}) we record this observation and its extensions to operators with small-range type, including the corresponding operator ideals and algebraicity modulo these ideals. These results supplement the main countable-type analysis.

\section{Preliminaries and Banach factorization}
Throughout, $\K$ is a complete nontrivially valued non-Archimedean field, $E$ is a Hausdorff locally convex space over $\K$, and $\Lop(E,F)$ denotes the continuous linear maps between locally convex spaces; $\Lop(E)=\Lop(E,E)$. By a seminorm we mean a non-Archimedean seminorm. For background we refer to van Rooij \cite{vanRooij} and P\'erez-Garc\'ia and Schikhof \cite{PGS}.

A set $C\subset E$ is \emph{compactoid} if, for every zero-neighbourhood $U$, there is a finite set $H\subset E$ with $C\subset U+\co(H)$, where $\co(H)$ is the set of all finite linear combinations of elements of $H$ with coefficients of absolute value at most one. A continuous linear map is called compactoid if it maps some zero-neighbourhood into a compactoid set. This is the terminology of \cite[Definition 8.1.1]{PGS}. Following \cite[Section 8.8, p. 334]{PGS}, an operator is called \emph{semicompact} if it maps some zero-neighbourhood into a compactoid Banach disk, and \emph{compact} if it maps some zero-neighbourhood into a complete compactoid set. A Banach disk $D$ is a bounded absolutely convex set for which $(\spn D,p_D)$, with its Minkowski norm, is Banach (\cite[p. 414]{PGS}). Compact operators are semicompact, and semicompact operators are compactoid. These three notions coincide for quasicomplete targets, in particular for Banach and Fr\'echet targets. We use ``relatively compact'' exclusively in the ordinary topological sense. Compactoidity is stronger than merely mapping bounded sets into compactoids (see \cite[Definition 8.3.1 and Remark 8.3.3]{PGS}).

A compactoid set is bounded \cite[Theorem 3.8.4(vii)]{PGS} and its closure is compactoid. For the latter, choose an open absolutely convex zero-neighbourhood $U$ contained in the given neighbourhood. The sum $U+\co(H)$ is an open additive subgroup, hence closed. Continuous linear maps preserve compactoids \cite[Theorem 3.8.4(iv)]{PGS}. Compactoid operators form an operator ideal by \cite[Theorem 8.1.3(i), (v)]{PGS}. For Banach spaces they are precisely the operator-norm limits of finite-rank operators, by \cite[Theorem 4.39, p. 142]{vanRooij}; equivalent characterizations, including a convergent series of rank-one operators, are given in \cite[Theorem 4.40, pp. 143--145]{vanRooij}. In these results van Rooij uses the term compact for what we call compactoid. If $\K$ is locally compact, compactoid subsets of a Banach space are relatively compact (cf. \cite[Theorem 3.8.3]{PGS}): the absolutely convex hull of a finite set is compact, so the defining inclusion implies total boundedness; completeness then gives compact closure. These observations will also be used for nonmetrizable spaces.

For $T\in\Lop(E)$, write $\{T\}'=\{A\in\Lop(E):AT=TA\}$. A closed linear subspace $M$ is hyperinvariant for $T$ if $A(M)\subset M$ for every $A\in\{T\}'$. It is nontrivial if $0\ne M\ne E$. In particular, for $x\ne0$, the closure of $\{T\}'x$ is a nonzero hyperinvariant subspace. The closure of $\{P(T)x:P\in\K[X]\}$ is invariant for $T$.

\begin{lemma}\label{lem:factor}
Suppose that $T\in\Lop(E)$ is compactoid and $E$ is sequentially complete, or that $T$ is semicompact without any completeness assumption on $E$. There exist a Banach space $F$ and continuous linear maps $J:E\to F$ and $S:F\to E$ such that
\begin{equation}\label{eq:factor}
 T=SJ,\qquad \cl{J(E)}=F,\qquad S\text{ is injective and compactoid}.
\end{equation}
If $T$ is injective, then $J$ is injective. If $T(E)$ is dense in $E$, then $V=JS$ has dense range in $F$. In all cases $V$ is compactoid and $VJ=JT$.
\end{lemma}
\begin{proof}
In the semicompact case choose a zero-neighbourhood $U$ and a compactoid Banach disk $D$ with $T(U)\subset D$. The map $J_0:E\to (\spn D,p_D)$, $J_0x=Tx$, is continuous, since $p_D(Tx)\le1$ on $U$. Let $F$ be the closure of $J_0(E)$ in this Banach space, let $J:E\to F$ be $J_0$ with restricted codomain, and let $S:F\to E$ be inclusion. The inclusion is continuous by boundedness of $D$, and compactoid: the closed unit ball of $(\spn D,p_D)$ lies in $\lambda D$ for any fixed $\lambda\in\K$ with $|\lambda|>1$, directly by the gauge definition. Hence $T=SJ$, $J(E)$ is dense, and $S$ is injective and compactoid. The other assertions follow from these properties and $VJ=JT$.

For a sequentially complete $E$, choose a continuous seminorm $p$ for which $T(U_p)$ is compactoid, where $U_p=\{x:p(x)\le1\}$. For each continuous seminorm $q$ there is a finite $C_q$ with
\begin{equation}\label{eq:bound}
 q(Tx)\le C_q p(x)\qquad(x\in E).
\end{equation}
Indeed, boundedness gives a uniform bound on $U_p$. Choose $\pi\in\K$ with $0<|\pi|<1$. For $p(x)>0$ scale $x$ by a suitable power of $\pi$ to bring its $p$-seminorm into $(|\pi|,1]$; the loss is at most $|\pi|^{-1}$. For $p(x)=0$, arbitrary scalings show that $q(Tx)=0$. In particular $\ker p\subset\ker T$.

Let $B$ be the Banach completion of $E/\ker p$ with its induced norm, and let $j:E\to B$ be the canonical map. By \eqref{eq:bound}, $T$ factors through $j(E)$ and extends uniquely to a continuous linear map $S_0:B\to E$. Indeed, approximate each $b\in B$ by a sequence from $j(E)$; its images form a Cauchy sequence in every continuous seminorm of $E$. Sequential completeness and the Hausdorff property give a well-defined limit, and \eqref{eq:bound} passes to the limit. Thus $S_0j=T$.

If $b$ belongs to the closed unit ball of $B$, the approximating sequence can be chosen in $j(U_p)$: any approximation with error less than one has norm at most one. Consequently
\[
 S_0(B_B)\subset\cl{T(U_p)}^{\,E},
\]
where $B_B$ denotes the closed unit ball of $B$. Hence $S_0$ is compactoid.

Put $N=\ker S_0$, $F=B/N$, let $Q:B\to F$ be the quotient map, and define $S(Qb)=S_0b$, $J=Qj$. Since $N$ is closed, $F$ is Banach \cite[Theorem 2.1.6(i)]{PGS}. The map $S$ is continuous and injective. It is compactoid because each member of the unit ball of $F$ has a representative of norm at most two in $B$, and $S_0$ maps every fixed-radius ball into a compactoid set. The remaining claims follow from $SJ=T$, from the density of $j(E)$ in $B$, and from $VJ=JT$. In particular, if $T(E)$ is dense then $V(F)$ contains the dense set $JT(E)$.
\end{proof}

\begin{remark}\label{rem:incomplete}
Sequential completeness can be dropped in Lemma~\ref{lem:factor} if $T$ maps a zero-neighbourhood into a relatively compact set (with compact closure in $E$). In that case, for each $b\in B$, an approximating sequence may be taken in a fixed-radius ball. Its image under $T$ is Cauchy and lies in a compact subset of $E$, so it converges in $E$. The rest of the proof is unchanged. Over a locally compact field, the relatively compact image is also compactoid. In a general, possibly incomplete locally convex space, we do not identify relatively compact and compactoid operators. In addition, for a sequentially complete Hausdorff target $F$, every compactoid map into $F$ is semicompact: the closure of the absolutely convex hull of its compactoid zero-neighbourhood image is bounded, compactoid and sequentially complete, and therefore a Banach disk by \cite[p. 414]{PGS}. The converse holds without completeness.

For a nonzero linear operator into a Hausdorff locally convex space, mapping a zero-neighbourhood into a relatively compact set actually forces $\K$ to be locally compact. Indeed, choose $x$ with $Tx\ne0$ and scale $x$ so that $\{\lambda x:|\lambda|\le1\}$ lies in the neighbourhood. Its image is homeomorphic to the closed unit ball of $\K$ and is closed in the target (use the closedness of the one-dimensional subspace and its usual topology). Relative compactness implies compactness of this ball. Thus the topological condition cannot by itself extend the nonzero case to non-locally compact fields.
\end{remark}

\section{Transfer of invariant and hyperinvariant subspaces}
For clarity, a normed space is of countable type if it has a countable subset with dense linear span. For an arbitrary locally convex space the monograph uses the term when each continuous seminorm gives a normed quotient of countable type \cite[Definition 4.2.1]{PGS}. Having a countable subset with dense linear span (strict countable type, \cite[Definition 4.2.3, p. 176]{PGS}) always implies countable type: the image of that subset has dense linear span in each normed quotient. The converse holds for metrizable spaces but not in general.

For later use, recall why the range of every compactoid operator $R:E\to F$ is of countable type. Choose a zero-neighbourhood $U$ for which $R(U)$ is compactoid. The set $U$ is absorbing, independently of completeness (\cite[Theorem 3.3.15(i)]{PGS}), so $R(E)=\spn R(U)$. Theorem 4.2.2 of \cite{PGS} gives countable type of this span with its induced topology; this is precisely Theorem 8.1.3(iv) there. Moreover, $R(U)$ is bounded by Theorem 3.8.4(vii), so every compactoid operator maps a zero-neighbourhood into a bounded set.

The following elementary lemma is a special case of \cite[Theorems 8.1.3(iv) and 4.2.13(ii)]{PGS}; we include its short normed-space proof.

\begin{lemma}\label{lem:countable}
If a compactoid operator $V$ on a Banach space $F$ has dense range, then $F$ is of countable type.
\end{lemma}
\begin{proof}
For each $n\ge1$ choose a finite set $H_n\subset F$ such that $V(B_F)\subset\{y:\|y\|\le2^{-n}\}+\co(H_n)$. Then $V(B_F)$ is contained in the closed linear span of $\bigcup_n H_n$. Scaling shows the same for $V(F)$, whose closure is $F$.
\end{proof}

The hyperinvariant part of the following lemma adapts Hoover's quasi-similarity argument \cite[Theorem 2.1]{Hoover} to Hausdorff locally convex spaces, with only density of $J(E)$ and injectivity of $S$ required. For ordinary invariance we use the product identities $T=SJ$ and $V=JS$.

\begin{transferlemma}\label{lem:transfer}
Let $E,F$ be Hausdorff locally convex spaces over $\K$, let $J:E\to F$ and $S:F\to E$ be continuous linear maps with $J(E)$ dense in $F$ and $S$ injective, and put $T=SJ$ and $V=JS$. If $V$ has a nontrivial closed hyperinvariant subspace, then so does $T$. The analogous assertion holds for invariant subspaces.
\end{transferlemma}
\begin{proof}
Let $M$ be such a hyperinvariant subspace of $F$, and choose $0\ne y\in M$. Then $x=Sy\ne0$. For each $A\in\{T\}'$, the continuous operator $W_A=JAS$ commutes with $V$, since
\[
 W_AV=JATS=JTAS=VW_A.
\]
Hence $JAx=W_Ay\in M$ for every $A\in\{T\}'$. Since $\{T\}'$ is a unital algebra, its orbit $\{T\}'x$ is a linear subspace, contains $x$, and is invariant under every member of $\{T\}'$. Its closure $H=\cl{\{T\}'x}^{\,E}$ is therefore nonzero and hyperinvariant. The orbit lies in the closed set $J^{-1}(M)$, hence so does $H$. Density of $J(E)$ and closedness and properness of $M$ imply that $J^{-1}(M)\ne E$. Thus $H$ is nontrivial. For ordinary invariance, $J^{-1}(M)$ itself is invariant by $VJ=JT$; it contains the nonzero vector $Sy$, since $JSy=Vy\in M$.
\end{proof}

\begin{theorem}\label{thm:equiv}
Fix a complete nontrivially valued non-Archimedean field $\K$. The following assertions are equivalent:
\begin{enumerate}
\item[(i)] Every nonzero compactoid operator on $\cnull$ has a nontrivial closed hyperinvariant subspace.
\item[(ii)] Every nonzero compactoid operator on any infinite-dimensional Banach space over $\K$ has a nontrivial closed hyperinvariant subspace.
\item[(iii)] Every nonzero compactoid operator on any infinite-dimensional Fr\'echet space over $\K$ has a nontrivial closed hyperinvariant subspace.
\item[(iv)] Every nonzero compactoid operator on any infinite-dimensional sequentially complete Hausdorff locally convex space over $\K$ has a nontrivial closed hyperinvariant subspace.
\item[(v)] Every nonzero semicompact operator on any infinite-dimensional Hausdorff locally convex space over $\K$ has a nontrivial closed hyperinvariant subspace.
\item[(vi)] Every nonzero compact operator on any infinite-dimensional Hausdorff locally convex space over $\K$ has a nontrivial closed hyperinvariant subspace.
\end{enumerate}
The equivalence also holds when ``hyperinvariant'' is replaced by ``invariant'' throughout.
\end{theorem}
\begin{proof}
It suffices to prove (i)$\Rightarrow$(iv) and (i)$\Rightarrow$(v): the other implications follow from the inclusions of operator classes and their coincidence for Banach targets. Let $E$ and $T$ be as in (iv) or (v). If $T$ is not injective, $\ker T$ is a nontrivial closed hyperinvariant subspace. If its range is not dense, $\cl{T(E)}$ is such a subspace. We may therefore assume $T$ is injective with dense range.

Use the maps $J,S$ of Lemma~\ref{lem:factor}. The Banach space $F$ is infinite-dimensional because $J$ is injective. The operator $V=JS$ is injective, compactoid and has dense range. By Lemma~\ref{lem:countable}, $F$ is of countable type. Every infinite-dimensional non-Archimedean Banach space of countable type is topologically linearly isomorphic to $\cnull$ \cite[Corollary 2.3.9]{PGS}. Thus (i) yields a nontrivial closed hyperinvariant subspace $M$ for $V$ in $F$.

Lemma~\ref{lem:transfer} gives the required hyperinvariant subspace of $T$. Notice that we have not asserted that an arbitrary member of the commutant of $T$ extends to the Banach completion of $E$. For ordinary invariance the same reduction and the invariant-subspace version of the transfer lemma apply.
\end{proof}

\begin{remark}
The reduction in Theorem~\ref{thm:equiv} relies on countable type of the auxiliary Banach space. In fact, in the dense-range case the original space is also of countable type in the sense of Definition 4.2.1 of \cite{PGS}, by Theorems 8.1.3(iv) and 4.2.13(ii) there. No assumption of strict countable type or metrizability of the original space is needed. Also, a nonzero compactoid operator on a locally convex space with no continuous norm cannot be injective, by \eqref{eq:bound}, and hence always has a nontrivial closed hyperinvariant subspace.
\end{remark}

\section{Fredholm theory and the positive-radius case}
For a bounded operator $V$ on a Banach space, let $\rad(V)=\lim_n\|V^n\|^{1/n}$. The limit exists by submultiplicativity. This is the norm-growth radius; it need not agree with the supremum of the absolute values of the spectrum over $\K$.

\begin{proposition}\label{prop:polynomial}
Let $E$ be Hausdorff locally convex, let $T\in\Lop(E)$ satisfy the hypotheses of Lemma~\ref{lem:factor}, and let $P\in\K[X]$ with $P(0)\ne0$. Then $P(T)$ has closed range and
\[
 \dim\ker P(T)=\dim(E/P(T)E)<\infty.
\]
In particular, if $P(T)$ is not bijective, its kernel is nonzero. If $E$ is infinite-dimensional, any such nonzero kernel is a nontrivial closed hyperinvariant subspace of $T$.
\end{proposition}
\begin{proof}
Write $P(T)/P(0)=I+C$, where $C=TQ(T)$ for a polynomial $Q$. If $T=SJ$ is a factorization as in Lemma~\ref{lem:factor}, then $C=S(JQ(T))$ factors through the same Banach space with a compactoid second factor. The next proposition and its proof apply to this factorization, yielding closed range and index zero. The kernel is closed by continuity and hyperinvariant because $P(T)$ commutes with every member of $\{T\}'$.
\end{proof}

\begin{remark}
The proposition does not infer the existence of a nonzero spectral value from index zero. The latter requires a separate spectral argument, supplied below by Serre's entire Fredholm determinant and resolvent. For sequentially complete spaces, the algebraic index-zero result (without closedness asserted) is exactly \cite[Corollary 3, p. 116]{SliwaFred}; for stability of Fredholm operators between Fr\'echet spaces see \cite[Theorem 4]{SliwaFred}. The factorization proof below is consistent with the classical locally convex Riesz theory in \cite[Section 5.2]{Gruson}.
\end{remark}

\begin{proposition}\label{prop:closed}
Let $X$ be Hausdorff locally convex. Suppose $C=SJ$, where $F$ is Banach, $J:X\to F$ is continuous linear and $S:F\to X$ is continuous linear and compactoid. Then $I_X+C$ is Fredholm of index zero, including closedness of its range. In particular this holds if $C$ satisfies the hypotheses of Lemma~\ref{lem:factor}.
\end{proposition}
\begin{proof}
Put $V=JS$. This is compactoid by the ideal property. By \cite[Corollary 3]{SliwaFred}, the Banach-space operator $I+V$ has index zero. Its range is closed: if $D$ is a finite-dimensional algebraic complement, the continuous surjection $F\times D\to F$, $(y,d)\mapsto (I+V)y+d$, is open by the Banach open mapping theorem \cite[Theorem 2.1.17]{PGS}. The inverse image of $(I+V)F$ is $F\times\{0\}$, hence $(I+V)F$ is closed. We have the exact range identity
\[
 (I_X+C)X=J^{-1}\big((I_F+V)F\big).
\]
The inclusion from left to right follows from $J(I+C)=(I+V)J$. Conversely, if $Jx=(I+V)y$, then $(I+C)(x-Sy)=x$, since $(I+C)S=S(I+V)$. The range identity proves closedness. The intertwining relations $J(I+C)=(I+V)J$ and $S(I+V)=(I+C)S$ show that $J$ and $-S$ induce maps between the kernels in opposite directions, and between the algebraic cokernels in opposite directions. On the kernels, $-SJ=-C$ and $-JS=-V$ act as the respective identities. On the cokernels they also act as the identities, since $x+Cx$ and $y+Vy$ belong to the corresponding ranges. Thus
\[
 \ker(I+C)\cong\ker(I+V),\qquad
 X/(I+C)X\cong F/(I+V)F.
\]
The algebraic index-zero result for the Banach operator transfers to $X$.
\end{proof}

\begin{lemma}\label{lem:fredholm}
Let $G=\cnull$ and let $V$ be compactoid on $G$. If $\rad(V)>0$, there exists an irreducible polynomial $P\in\K[X]$, $P(0)\ne0$, with $0<\dim\ker P(V)<\infty$.
\end{lemma}
\begin{proof}
Let $P_n$ be the coordinate projection onto the first $n$ coordinates. By the compactoid characterization in \cite[Lemma 8.1.4 and Theorem 8.1.9(ii)]{PGS}, or directly as follows, compactoidity implies $\|(I-P_n)V\|\to0$: approximate $V(B_G)$ by the absolutely convex hull of a finite set up to an arbitrarily small ball, and use uniform convergence of the coordinate tails on that finite set. Thus $V$ is a norm limit of finite-rank operators. In Serre's terminology it is completely continuous. Since $G$ already has its canonical orthonormal basis, the Fredholm theory on $c(I)$ in \cite[Propositions 7, 10 and 11]{Serre} applies without any discreteness assumption. Indeed, the determinant and resolvent coefficient estimates on pp. 76 and 79 use the bound $r_1\cdots r_m$, where $(r_j)$ is the nonincreasing arrangement of the matrix row bounds and $r_j\to0$. For each $s>0$ one has $s^m r_1\cdots r_m\to0$. These estimates and the determinant identity in Proposition 11 use only completeness and the non-Archimedean absolute value; they do not require the existence of a uniformizer.

Write $D(z)=\det(I-zV)$ for the entire Fredholm determinant. The formal Fredholm resolvent $D(z)(I-zV)^{-1}$ is entire with values in $\Lop(G)$. If $D\equiv1$, its coefficients are $V^n$, so $s^n\|V^n\|\to0$ for every real $s>0$. This would imply $\rad(V)=0$. Therefore $D$ is nonconstant.

Write $D(z)=\sum_{n\ge0}d_nz^n$ with $d_0=1$. Choose $r\in|\K^\times|$ so large that $|d_m|r^m>1$ for some $m\ge1$ with $d_m\ne0$. In the completion of an algebraic closure of $\K$, consider the closed disk of radius $r$ about zero. The largest index $j$ at which $\max_{n\ge0}|d_n|r^n$ is attained is positive. By the Weierstrass preparation theorem \cite[Theorem 3.5, p. 15]{Benedetto}, $D$ factors on this disk as a monic polynomial in $\K[z]$ of degree $j$ times a zero-free analytic factor. The zero-counting statement \cite[Theorem 3.9, p. 17]{Benedetto} gives a zero $\alpha$ of $D$ algebraic over $\K$. These theorems apply to arbitrary complete non-Archimedean fields, in particular to discretely valued ones. Since $D(0)=1$, $\alpha\ne0$.

Over the finite complete extension $L=\K(\alpha)$ (cf. \cite[Theorems 1.2.10 and 2.1.14]{PGS}), extend $V$ by its matrix to $c_0(\N,L)$. The coordinate-tail norm estimates are unchanged; the Fredholm determinant is the same by its matrix construction. By \cite[Proposition 11, p. 80]{Serre}, $D(\alpha)=0$ implies that $I-\alpha V_L$ is not invertible. Proposition~\ref{prop:closed} applied over $L$ shows that it has closed range and index zero. If its kernel were zero, it would be bijective and hence invertible by the Banach open mapping theorem. Consequently there is a nonzero vector $v$ satisfying $V_Lv=\alpha^{-1}v$. Let $P$ be the minimal polynomial of $\alpha^{-1}$ over $\K$. Then $P(0)\ne0$ and $P(V_L)v=0$. Expanding $v$ in a finite $\K$-basis of $L$ gives vectors in $G$ (the coordinate functionals of that basis are continuous), at least one nonzero, annihilated by $P(V)$. Hence $\ker P(V)\ne0$.

Proposition~\ref{prop:polynomial} gives $\dim\ker P(V)<\infty$. This proves the lemma, in arbitrary characteristic. It is an alternative proof of the polynomial-kernel conclusion in \cite[Theorem 3.5]{EKB} for countable type.
\end{proof}

\begin{theorem}\label{thm:banachradius}
Let $G$ be any Banach space over $\K$ and let $V\in\Lop(G)$ be compactoid with $\rad(V)>0$. There exists an irreducible $P\in\K[X]$ with $P(0)\ne0$ and
\[
 0<\dim\ker P(V)<\infty.
\]
In particular, if $G$ is infinite-dimensional, this kernel is a nontrivial closed hyperinvariant subspace of $V$.
\end{theorem}
\begin{proof}
Put $Y=\cl{V(G)}$. By \cite[Theorems 8.1.3(iv) and 4.2.13(ii)]{PGS}, $Y$ is a Banach space of countable type. It is invariant for $V$, and the restriction $A=V|_Y:Y\to Y$ is compactoid by \cite[Theorem 8.1.3(ii), (iii)]{PGS}. This choice of $Y$ avoids any assumption that $G$ has a basis or that $Y$ is complemented.

Since $V(G)\subset Y$, for $n\ge1$,
\[
 \|A^n\|\le\|V^n\|,\qquad
 \|V^{n+1}\|\le\|A^n\|\,\|V\|.
\]
Taking roots gives $\rad(A)=\rad(V)>0$. If $Y$ is infinite-dimensional, it is isomorphic to $\cnull$ and Lemma~\ref{lem:fredholm} gives a polynomial $P$ with the required kernel for $A$. If $Y$ is finite-dimensional, write the characteristic polynomial of $A$ as $X^mP(X)$ with $P(0)\ne0$. The positive radius excludes nilpotence. If $P(A)$ were invertible, the Cayley--Hamilton theorem would imply $A^m=0$, a contradiction. Thus $P(A)$ is not invertible. Factoring $P$ into irreducibles, some irreducible factor $R$ has $R(A)$ noninvertible (otherwise their product would be invertible). Since $Y$ is finite-dimensional, $\ker R(A)\ne0$. Replace $P$ by this factor; it still has $P(0)\ne0$.

Finally, $P(0)\ne0$ implies $\ker P(V)\subset V(G)\subset Y$, since $P(V)x=0$ expresses $x$ as a linear combination of positive powers of $V$ applied to $x$. Therefore $\ker P(V)=\ker P(A)$. It is closed and invariant under the entire commutant because $P(V)$ commutes with every member of $\{V\}'$.
\end{proof}

\begin{remark}
Theorem~\ref{thm:banachradius} is a basis-free version of \cite[Theorem 3.5 and Corollary 3.6]{EKB}. It is a consequence of classical Fredholm theory; in particular, the countable-type reduction is closely related to \cite[Section 4.3]{Gruson}. We include the proof to make the argument for locally convex spaces explicit. The locally convex Riesz theory for semi-compact endomorphisms was already discussed in \cite[Section 5.2, p. 93]{Gruson}; the contribution here is the combination of Banach factorization with the full-commutant transfer argument (cf. Hoover \cite[Theorem 2.1]{Hoover}) in the non-Archimedean locally convex setting, rather than the existence of nonzero algebraic spectral values itself.
\end{remark}

\subsection*{A radius for compactoid endomorphisms of locally convex spaces}
For $T$ as in Lemma~\ref{lem:factor}, put $\rho(T)=\rad(JS)$ for a Banach factorization $T=SJ$. In fact, the definition is independent of the factorization. For two such factorizations $T=S_iJ_i$ through Banach spaces $F_i$, let $V_i=J_iS_i$, $A=J_2S_1:F_1\to F_2$ and $B=J_1S_2:F_2\to F_1$. Then $BA=V_1^2$ and $AB=V_2^2$. The identity $\rad(AB)=\rad(BA)$ follows from $(AB)^{n+1}=A(BA)^nB$ and the reverse estimate. Thus $\rad(V_1)=\rad(V_2)$. This argument works for any continuous Banach factorizations, and in the Banach case $\rho(T)=\rad(T)$ by choosing the identity as one factor.

\begin{theorem}\label{thm:lcsradius}
Let $E$ be an infinite-dimensional Hausdorff locally convex space and $T\in\Lop(E)$ satisfy the hypotheses of Lemma~\ref{lem:factor}. If $\rho(T)>0$, then there exists an irreducible $P\in\K[X]$, $P(0)\ne0$, such that $0<\dim\ker P(T)<\infty$. If $\rho(T)=0$, then for every continuous seminorm $q$, every $x\in E$ and every real $s>0$,
\[
 s^n q(T^n x)\longrightarrow0.
\]
\end{theorem}
\begin{proof}
Take the factorization $T=SJ$ and set $V=JS$. If $\rho(T)>0$, apply Theorem~\ref{thm:banachradius} to $V$. The intertwining relations are
\[
 TS=SV,\qquad VJ=JT.
\]
In particular $S(\ker P(V))\subset\ker P(T)$ and $J(\ker P(T))\subset\ker P(V)$. The restriction of $S$ to $\ker P(V)$ is injective: if $Sy=0$ and $P(V)y=0$, then $Vy=0$ and $P(0)y=0$. Similarly, $J$ is injective on $\ker P(T)$. Thus the two kernels have the same finite positive dimension. The kernel in $E$ is closed and hyperinvariant.

If $\rho(T)=0$, use $T^n=SV^{n-1}J$ for $n\ge1$. For every $q$ there is a constant $C_q$ with $q(Sy)\le C_q\|y\|$. Hence $q(T^n x)\le C_q\|V^{n-1}\|\|Jx\|$, giving the asserted decay.
\end{proof}

\section{Locally compact fields and quasi-null vectors}
The following theorem extends the quasi-null-vector argument of \cite[Theorem 9]{BE} and \cite[Theorem 3.14]{EKB} to sequentially complete locally convex spaces. The additional composition in the proof ensures a relatively compact image in the original topology, without assuming completeness for general Cauchy nets.

\begin{theorem}\label{thm:quasinull}
Let $\K$ be locally compact and $E$ be an infinite-dimensional sequentially complete Hausdorff locally convex space over $\K$. Let $T\in\Lop(E)$ be nonzero and commute with a compactoid operator $C\in\Lop(E)$. Suppose there is a continuous seminorm $p$ such that $C(U_p)$ is compactoid, where $U_p=\{x:p(x)\le1\}$, and a vector $x_0\in E$ such that
\[
 Cx_0\ne0,\qquad \liminf_{n\to\infty}p(T^nx_0)^{1/n}=0.
\]
Then $T$ has a nontrivial closed hyperinvariant subspace.
\end{theorem}
\begin{proof}
Suppose that no such subspace exists. Then $T$ is injective and, writing $\mathcal A=\{T\}'$, the subspace $\mathcal A x$ is dense for each $x\ne0$.

Since $Cx_0\ne0$ and $\ker C$ is closed and proper, choose $A_0\in\mathcal A$ such that $CA_0Cx_0\ne0$. Put $D=CA_0C$. Then $D$ commutes with $T$ and $Dx_0\ne0$. By the proof of Lemma~\ref{lem:factor} before taking the quotient by $\ker S_0$, there are a Banach space $B$, a continuous map $j:E\to B$ with $\|jx\|=p(x)$, and a compactoid continuous map $S_0:B\to E$ such that $C=S_0j$. Thus
\[
 D=S_0(jA_0S_0)j.
\]
The middle factor $jA_0S_0$ is a compactoid endomorphism of $B$. Local compactness of $\K$ and completeness of $B$ imply that the closure of its unit-ball image is compact in $B$ \cite[Theorem 3.8.3]{PGS}. Its continuous image under $S_0$ is compact in $E$ and contains $D(U_p)$. In particular $D(U_p)$ is relatively compact in $E$. This reasoning avoids assuming that every compactoid set in a merely sequentially complete space has compact closure.

As $TDx_0\ne0$, choose a continuous seminorm $q$ with $q(TDx_0)>0$. Boundedness of $TD(U_p)$ gives a constant $L>0$ with $q(TDx)\le Lp(x)$ for every $x\in E$ (by the scaling argument in Lemma~\ref{lem:factor}). Also $p(x_0)>0$, because $Cx_0\ne0$. Choose $0<\varepsilon<p(x_0)$ with $L\varepsilon<q(TDx_0)$ and put
\[
 U=\{x\in E:p(x-x_0)\le\varepsilon\}.
\]
The set $U$ is open and closed. For $x\in U$, ultrametricity gives $p(x)=p(x_0)$ and $q(TDx)=q(TDx_0)$. Moreover $U$ is contained in a scalar multiple of $U_p$, so $Z=\cl{TD(U)}$ is compact in $E$ and excludes zero.

The sets $A^{-1}(U)$, $A\in\mathcal A$, cover $Z$. A finite subcover is given by $A_1,\ldots,A_m$ in $\mathcal A$. Starting at $u_0=x_0$, choose recursively
\[
 u_n=A_{i_n}TDu_{n-1}\in U.
\]
Each $A_i$ and $D$ commutes with $T$ (they need not commute with one another), so
\[
 u_n=A_{i_n}D\cdots A_{i_1}DT^n x_0.
\]
Since $A_iD(U_p)$ is bounded, a common $M\ge1$ satisfies $p(A_iDx)\le Mp(x)$ for all $x\in E$ and $1\le i\le m$. It follows that $p(u_n)\le M^n p(T^nx_0)$, whose liminf is zero. But $p(u_n)=p(x_0)>0$ for all $n$. This contradiction proves the theorem.
\end{proof}

\begin{remark}
For Banach spaces take $p=\|\cdot\|$ in Theorem~\ref{thm:quasinull}. This gives the nonzero version of \cite[Theorem 3.14]{EKB}. The exclusion of the zero operator is essential for hyperinvariance in spaces such as $\cnull$, where continuous rank-one operators act transitively on nonzero vectors. Both proofs use a finite covering of a compact closure; in our locally convex setting the additional composition $D=CA_0C$ ensures the relative compactness needed for this step.
\end{remark}

\begin{corollary}\label{cor:local}
Let $\K$ be locally compact. Every nonzero compactoid operator on an infinite-dimensional sequentially complete Hausdorff locally convex space over $\K$ has a nontrivial closed hyperinvariant subspace. The same conclusion holds for every nonzero semicompact (in particular, compact) operator, and for a nonzero operator mapping a zero-neighbourhood into a relatively compact set, on an arbitrary infinite-dimensional Hausdorff locally convex space over $\K$.
\end{corollary}
\begin{proof}
For the first assertion let $T$ be such an operator. If $\rho(T)>0$, apply Theorem~\ref{thm:lcsradius}. If $\rho(T)=0$, choose a continuous seminorm $p$ for which $T(U_p)$ is compactoid and a vector $x_0$ with $Tx_0\ne0$. The decay assertion of Theorem~\ref{thm:lcsradius} implies $\lim_n p(T^nx_0)^{1/n}=0$. Apply Theorem~\ref{thm:quasinull} with $C=T$.

For the second assertion, use the semicompact factorization of Lemma~\ref{lem:factor} or the relatively compact factorization in Remark~\ref{rem:incomplete} and the proof of Theorem~\ref{thm:equiv}, which requires no further completeness of the original space once that factorization is obtained. The first assertion applies to the resulting Banach space.
\end{proof}

\begin{corollary}\label{cor:lomo}
Let $E$ be infinite-dimensional and Hausdorff locally convex over a locally compact $\K$. If $T\in\Lop(E)$ commutes with a nonzero compactoid operator $C$ and $E$ is sequentially complete, or if it commutes with a nonzero semicompact operator $C$ without a completeness assumption on $E$, then $T$ has a nontrivial closed invariant subspace.
\end{corollary}
\begin{proof}
Apply Corollary~\ref{cor:local} to $C$. A nontrivial closed hyperinvariant subspace for $C$ is invariant for $T$ since $T\in\{C\}'$.
\end{proof}

\begin{corollary}\label{cor:positivecomm}
Let $\K$ be arbitrary complete and nontrivially valued, and let $E$ be infinite-dimensional and Hausdorff locally convex. If $T\in\Lop(E)$ commutes with an operator $C$ satisfying the hypotheses of Lemma~\ref{lem:factor} and $\rho(C)>0$, then $T$ has a nontrivial finite-dimensional closed invariant subspace.
\end{corollary}
\begin{proof}
Theorem~\ref{thm:lcsradius} applied to $C$ gives a nonzero finite-dimensional kernel $\ker P(C)$. Since $T$ commutes with $C$, it leaves this kernel invariant. No local compactness hypothesis is needed.
\end{proof}

\section{Algebraic obstructions and commutants}

We give an algebraic refinement of the preceding criterion. The next proposition extends the necessary condition in \cite[Theorem 3.10]{EKB} under a separating-dual assumption, and identifies the precise obstruction for algebraic endomorphisms. No completeness of the underlying locally convex space is needed in this proposition. For nonalgebraic operators the necessary condition of \cite[Theorem 3.10]{EKB} holds automatically, since no nonzero polynomial annihilates such an operator.

\begin{proposition}\label{prop:alg}
Let $E$ be an infinite-dimensional Hausdorff locally convex space whose continuous dual separates points. Let $T\in\Lop(E)$ be algebraic, and let $m_T\in\K[X]$ be its monic minimal polynomial. Then $T$ has a nontrivial closed hyperinvariant subspace if and only if $m_T$ is reducible in $\K[X]$.
\end{proposition}
\begin{proof}
If $m_T=ab$ with $a,b$ of positive degree, then $a(T)$ and $b(T)$ are nonzero by minimality, whereas $a(T)b(T)=0$. Thus $\ker a(T)$ is a nonzero proper closed hyperinvariant subspace. This direction does not require separation by the dual.

Suppose instead that $m_T$ is irreducible and let $L=\K[X]/(m_T)$, a finite field extension of $\K$. Equip it with its unique extended valuation. Write $\alpha$ for the class of $X$ and let $L$ act on $E$ by $a\cdot x=a(T)x$, where a polynomial representative of $a$ is used. This is well-defined, and for each fixed $a\in L$ the action on $E$ is continuous. The continuous $\K$-linear maps commuting with $T$ are exactly the continuous $L$-linear maps on $E$ (with its original topology).

Fix $0\ne x\in E$. Choose a continuous $\K$-linear functional $f:E\to\K$ with $f(x)\ne0$. Let $\ell:L\to\K$ be any nonzero $\K$-linear functional. The pairing $(a,b)\mapsto\ell(ab)$ is nondegenerate: if $b\ne0$, then $Lb=L$, so $\ell(Lb)\ne0$. Thus for every $y\in E$ there is a unique $F(y)\in L$ such that
\[
 \ell(a F(y))=f(a\cdot y)\quad (a\in L).
\]
Taking a finite $\K$-basis of $L$ in this equality shows that $F$ is continuous and $\K$-linear. For $c\in L$, the same equality gives $F(c\cdot y)=cF(y)$ by nondegeneracy. Hence $F$ is $L$-linear, and $F(x)\ne0$ because $\ell(F(x))=f(x)\ne0$.

For each $z\in E$, the map
\[
 A(y)=\big(F(y)/F(x)\big)\cdot z
\]
is continuous (expand the coefficient in a finite basis of $L$), is $L$-linear, and satisfies $Ax=z$. Therefore the commutant acts transitively on nonzero vectors and has no nontrivial common invariant subspace. This proves the assertion. The argument works also for inseparable extensions, since it uses a nonzero linear functional $\ell$, not the field trace.
\end{proof}

\begin{theorem}\label{thm:dichotomy}
Let $E$ be an infinite-dimensional Hausdorff locally convex space, and let $T\in\Lop(E)$ commute with an endomorphism $C$ satisfying the hypotheses of Lemma~\ref{lem:factor} and $\rho(C)>0$. Then either $T$ has a nontrivial closed hyperinvariant subspace or $T$ is annihilated by an irreducible polynomial in $\K[X]$.

More precisely, $T$ has no nontrivial closed hyperinvariant subspace if and only if the continuous dual of $E$ separates points and $T$ is annihilated by an irreducible polynomial in $\K[X]$.
\end{theorem}
\begin{proof}
By Theorem~\ref{thm:lcsradius}, there exists a nonzero finite-dimensional subspace $M=\ker P(C)$, with $P(0)\ne0$, which is hyperinvariant for $C$. In particular $T(M)\subset M$. Let $p$ be the characteristic polynomial of $T|_M$. The Cayley--Hamilton theorem gives $M\subset\ker p(T)$. If $p(T)\ne0$, this kernel is nonzero, proper, closed and hyperinvariant for $T$.

Otherwise $T$ is algebraic. If its minimal polynomial is reducible, the first part of the proof of Proposition~\ref{prop:alg} applies without any assumption on the dual. The remaining possibility is that the minimal polynomial is irreducible. Under the separating-dual hypothesis Proposition~\ref{prop:alg} excludes hyperinvariant subspaces in that case. Being annihilated by an irreducible polynomial is equivalent to having that polynomial (up to a scalar) as the minimal polynomial, since the minimal polynomial of an endomorphism of a nonzero space is not constant.

It remains to show that failure of the separating-dual hypothesis implies existence of a nontrivial closed hyperinvariant subspace. We use a more general observation, valid for every Hausdorff locally convex space $E$: the existence of any nonzero compactoid endomorphism $C$ implies $E'\ne\{0\}$. Indeed, the space $C(E)$ with its induced topology is of countable type by \cite[Theorem 8.1.3(iv)]{PGS}, and its dual separates points by \cite[Corollary 5.1.4]{PGS}. Choose $x\in E$ with $Cx\ne0$ and $f\in C(E)'$ with $f(Cx)\ne0$. Then $fC$ is a nonzero element of $E'$. This observation does not require positive radius or sequential completeness.

If $E'$ does not separate points, the subspace
\[
 N_E=\bigcap_{f\in E'}\ker f
\]
is therefore nonzero and proper. It is closed and invariant under every continuous linear endomorphism of $E$, because $fA\in E'$ whenever $f\in E'$ and $A\in\Lop(E)$. In particular it is hyperinvariant for $T$. This completes the characterization.
\end{proof}

\begin{remark}
In particular, over a locally compact field the separating-dual hypothesis is automatic by \cite[Theorems 1.2.8, 1.2.13 and Corollary 5.1.4]{PGS}. Theorem~\ref{thm:dichotomy} gives a partial answer to the question at the end of \cite{EKB} under the additional hypothesis $\rho(C)>0$ on the commuting compactoid operator $C$, not on $T$. Here the question is understood for a nonzero commuting compactoid operator on an infinite-dimensional space. It is not settled when the only available nonzero compactoid operators in the commutant have radius zero. The two radii must not be confused.
\end{remark}

\begin{corollary}
Let $E$ be infinite-dimensional and Hausdorff locally convex over a locally compact field. Suppose $T\in\Lop(E)$ and there exists a nonzero operator $C$ satisfying the hypotheses of Lemma~\ref{lem:factor} in the bicommutant $\{T\}''=\{A\in\Lop(E):AB=BA\text{ for every }B\in\{T\}'\}$. Then $T$ has a nontrivial closed hyperinvariant subspace. In particular this holds if $P(T)\ne0$ is compactoid and $E$ is sequentially complete, or if $P(T)\ne0$ is semicompact, for some $P\in\K[X]$.
\end{corollary}
\begin{proof}
Every member of $\{T\}'$ commutes with $C$. A nontrivial closed hyperinvariant subspace of $C$ provided by Corollary~\ref{cor:local} is therefore hyperinvariant for $T$. Each polynomial in $T$ belongs to the bicommutant.
\end{proof}

\begin{corollary}\label{cor:semi}
Let $E$ be an infinite-dimensional sequentially complete Hausdorff locally convex space over a locally compact field, and suppose all bounded subsets of $E$ are compactoid. If $T\in\Lop(E)$ and $T-\lambda I$ maps some zero-neighbourhood into a bounded set for some $\lambda\in\K$, then $T$ has a nontrivial closed invariant subspace. If $T\ne\lambda I$, it has a nontrivial closed hyperinvariant subspace.
\end{corollary}
\begin{proof}
The operator $T-\lambda I$ is compactoid. If it is nonzero, apply Corollary~\ref{cor:local} and use $\{T\}'=\{T-\lambda I\}'$. If $T=\lambda I$, every closed linear subspace is invariant; a nonzero one-dimensional subspace is closed and proper.
\end{proof}

\begin{remark}
The last corollary applies to nuclear Fr\'echet spaces, in particular nuclear non-Archimedean K\"othe spaces, because their bounded sets are compactoid \cite[Theorem 8.4.5($\delta$) and Corollaries 8.5.2--8.5.3]{PGS}. It also clarifies the relation with the negative results of Kasprzak \cite[Theorems 3.8 and 4.1]{Kasprzak}, with the corrections in \cite{KasprzakCorr}. Those theorems allow any complete nontrivially valued non-Archimedean scalar field, including locally compact ones. On a nuclear K\"othe space covered by those results over a locally compact field, an operator $T$ with no nontrivial closed invariant subspace cannot have $T-\lambda I$ locally bounded for any $\lambda\in\K$, by the preceding corollary. Thus the assertion that every continuous endomorphism is a scalar plus a locally bounded operator fails on those spaces. For a concrete example, the matrix $b_n^k=(k+1)^n$ satisfies the assumptions of Kasprzak's Theorem 4.1: its successive-coordinate ratios are constant in $n$ and $b_n^{k+1}/b_n^k\to\infty$. This is a nuclear K\"othe space. Also, Corollary~\ref{cor:lomo} asserts invariance for $T$, not hyperinvariance for $T$.
\end{remark}

\section{Extension of the scalar field}
\begin{theorem}\label{thm:extension}
Let $K_0$ be a complete nontrivially valued non-Archimedean field and let $\K$ be any complete isometric valued extension. Let $T_0$ be a bounded operator on $c_0(\N,K_0)$ and let $T$ be its extension by the same matrix to $c_0(\N,\K)$. If $M_0$ is a nontrivial closed hyperinvariant subspace of $T_0$, then $M=\overline{\operatorname{span}_{\K}M_0}$ is a nontrivial closed hyperinvariant subspace of $T$. No basis or spherical completeness hypothesis is imposed on the field extension.
\end{theorem}
\begin{proof}
The extension is bounded because its matrix has uniformly bounded entries and each column tends to zero, as holds for the matrix of $T_0$. Its coordinate sums converge for every member of $c_0(\N,\K)$, and the resulting map has image in $c_0(\N,\K)$ by density of the finite sequences. The subspace $M$ is nonzero. It is proper: $c_0(\N,K_0)/M_0$ is a nonzero Banach space of countable type, so its continuous dual separates points \cite[Theorem 4.2.13(i) and Corollary 5.1.4]{PGS}. Composing a nonzero functional on this quotient with the quotient map gives a nonzero bounded $K_0$-linear functional on $c_0(\N,K_0)$ vanishing on $M_0$. By \cite[Theorem 2.5.11]{PGS}, its bounded sequence of coordinate coefficients defines a nonzero bounded $\K$-linear functional on $c_0(\N,\K)$ vanishing on $M$. Indeed, the extension is given by the convergent series of coordinate products; it annihilates the $\K$-linear span of $M_0$ and hence its closure.

Fix any bounded $\K$-linear operator $A$ commuting with $T$; if $A=0$ there is nothing to prove. Write $(a_{ij})_{i,j\in\N}$ for its matrix. Let
\[
 Y=\overline{\operatorname{span}_{K_0}\{a_{ij}:i,j\in\N\}}^{\,\K}.
\]
This is a nonzero Banach space of countable type over $K_0$, even if $\K$ itself is not of countable type and has no orthogonal basis over $K_0$. Since all the entries of $A$ lie in $Y$, continuity and approximation by finite sequences show that $A(c_0(\N,K_0))\subset c_0(\N,Y)$. Here $c_0(\N,Y)$ denotes the Banach space of sequences in $Y$ tending to zero in its norm. Also $T$ leaves $c_0(\N,Y)$ invariant, because its matrix entries belong to $K_0$ and $Y$ is a closed $K_0$-linear subspace of $\K$.

The countable-type isomorphism theorem \cite[Corollary 2.3.9 and Theorem 2.1.14]{PGS} gives a topological linear isomorphism $c_0(I,K_0)\to Y$ for a finite or countably infinite set $I$. Let $b_\eta\in Y$ and $\phi_\eta:Y\to K_0$ ($\eta\in I$) be its basis vectors and coordinate functionals. Both families are uniformly bounded in their respective norms, by boundedness of the isomorphism and its inverse. Write $\Phi_\eta:c_0(\N,Y)\to c_0(\N,K_0)$ for the map given by applying $\phi_\eta$ coordinatewise. For every $z\in c_0(\N,Y)$ the expansion
\[
 z=\sum_{\eta\in I} b_\eta\Phi_\eta z
\]
converges in the sup norm (as finite partial sums). Indeed the expansion converges on each coordinate, and the finite partial-sum maps on $Y$ have uniformly bounded norms; convergence on finitely many initial coordinates and the smallness of the remaining coordinates give uniform convergence.

The operator $A_\eta=\Phi_\eta A|_{c_0(\N,K_0)}$ is bounded and $K_0$-linear. Continuity of $\phi_\eta$ and the matrix expansion of $T$ give $\Phi_\eta Tz=T_0\Phi_\eta z$ for $z\in c_0(\N,Y)$. From $AT=TA$ it follows that $A_\eta T_0=T_0 A_\eta$. For $x\in M_0$, hyperinvariance gives $A_\eta x\in M_0$. The expansion of $Ax$ in terms of $b_\eta A_\eta x$ shows that $Ax\in M$. By $\K$-linearity and continuity, $A(M)\subset M$. As $A$ was arbitrary, $M$ is hyperinvariant for the full commutant of $T$.
\end{proof}

\begin{corollary}\label{cor:extension}
Suppose $K_0$ is locally compact and $\K$ is a complete isometric valued extension. Every nonzero compactoid operator on $c_0(\N,\K)$ whose matrix entries belong to $K_0$ has a nontrivial closed hyperinvariant subspace.
\end{corollary}
\begin{proof}
If the valuation induced on $K_0$ is trivial, compactoidity gives $\sup_{i>n,j}|t_{ij}|\to0$ by the coordinate-tail estimate in the proof of Lemma~\ref{lem:fredholm}. Since every nonzero entry has absolute value one, all rows beyond some index are zero. Thus $T$ has finite nonzero rank, and its range is a nontrivial closed hyperinvariant subspace. We may therefore assume the valuation on $K_0$ is nontrivial.

The matrix restricts to a nonzero bounded operator $T_0$ on $c_0(\N,K_0)$, since $K_0$ is complete and closed in $\K$. It is compactoid: the coordinate-tail norm estimates for the original operator give the same estimates for $T_0$, and its finite-row truncations converge to $T_0$ in norm (cf. \cite[Theorem 8.1.9(ii)]{PGS}). Over the locally compact field $K_0$, Corollary~\ref{cor:local} gives a nontrivial closed hyperinvariant subspace for $T_0$. Theorem~\ref{thm:extension} now gives the conclusion without any assumption on $\K$ as a Banach space over $K_0$.
\end{proof}

\begin{remark}
The countable-type argument in Theorem~\ref{thm:extension} concerns the matrix entries of one commuting operator at a time, not the extension field as a whole. The matrix assumption in Corollary~\ref{cor:extension} is substantive. For a general operator on a non-locally compact field it need not hold, even after choosing a locally compact subfield. In particular, $\mathbb C((t))$ with its $t$-adic valuation has residue characteristic zero and cannot contain a nontrivially valued locally compact subfield with the induced valuation. Thus the theorem does not resolve the general problem for discretely valued fields of infinite residue field. It does cover certain operators over non-locally compact extensions of $\mathbb Q_p$, such as operators with matrices over $\mathbb Q_p$ itself.
\end{remark}

\begin{corollary}\label{cor:separablefield}
Consider the assertion that every nonzero compactoid operator on $c_0(\N,K)$ has a nontrivial closed hyperinvariant subspace for every complete nontrivially valued non-Archimedean field $K$. This is equivalent to the same assertion restricted to fields separable in the valuation topology. The same holds for invariant subspaces, and separately for the class of discretely valued complete fields.
\end{corollary}
\begin{proof}
Given a matrix $(t_{ij})$ over $K$, choose $\pi\in K$ with $0<|\pi|<1$. Let $K_0$ be the closure in $K$ of the field generated by its prime field, $\pi$ and all the $t_{ij}$. This field is complete and nontrivially valued, with a countable dense subfield, since the generating set is countable. If $K$ is discretely valued, so is $K_0$. The same matrix defines a nonzero compactoid operator over $K_0$ by the coordinate-tail argument in Corollary~\ref{cor:extension}. Apply the assumed assertion over $K_0$ and then Theorem~\ref{thm:extension}. For ordinary invariance, extension of a closed invariant subspace by taking its closed $K$-linear span preserves invariance by continuity, and properness is proved by the same functional as in Theorem~\ref{thm:extension}.
\end{proof}

\section{Further questions}\label{sec:questions}
By Theorem~\ref{thm:equiv} and the spectral result, it remains to study injective compactoid operators on $\cnull$ with dense range and zero norm-growth radius. By Corollary~\ref{cor:separablefield}, for the problem over all fields it suffices to consider separable complete fields. On spaces not of countable type the compactoid case is immediate, as explained in Section~\ref{sec:cardinal}. In particular the following questions are equivalent, respectively, to their analogues on all sequentially complete Hausdorff locally convex spaces and to the semicompact and compact variants on arbitrary Hausdorff locally convex spaces.
\begin{problem}
Does every nonzero compactoid operator on $\cnull$ have a nontrivial closed invariant subspace when $\K$ is an arbitrary complete nontrivially valued non-Archimedean field?
\end{problem}
\begin{problem}
Does every such operator have a nontrivial closed hyperinvariant subspace? In particular, can the locally compact hypothesis in Corollary~\ref{cor:local} be replaced by spherical completeness or by discrete valuation alone?
\end{problem}
The finite-cover proof of Theorem~\ref{thm:quasinull} does not extend merely by replacing compact sets with compactoid sets: over a non-locally compact field a compactoid set need not be totally bounded. The scalar-extension result of Corollary~\ref{cor:extension} covers a restricted class of operators on some non-locally compact fields. For comparison, El Asri and Babahmed \cite{EB2023} obtained spectral quasinilpotent operators without nontrivial closed invariant subspaces on infinite-dimensional Banach spaces of countable type; compactoidity is not part of the statement of that result. No resolution of the general remaining cases is claimed here.

\section{Operators on spaces of uncountable type}\label{sec:cardinal}
As a complement to the preceding results, we consider spaces not of countable type. The countable-type assertion in \cite[Theorem 8.1.3(iv)]{PGS} already settles the compactoid problem on spaces not of countable type: the closure of the range of a nonzero compactoid endomorphism is a proper closed hyperinvariant subspace. We record a more general formulation which does not assume compactoidity.

Fix an infinite cardinal $\kappa$. A normed space is said to be \emph{of type at most $\kappa$} if it has a subset of cardinality at most $\kappa$ whose linear span is dense. A Hausdorff locally convex space is of type at most $\kappa$ if each normed quotient $E_p=E/\ker p$, for $p$ a continuous seminorm, has this property. At $\kappa=\aleph_0$ this is precisely \cite[Definition 4.2.1]{PGS}; in general it does not require a single dense linear span of cardinality at most $\kappa$ for the original topology.

\begin{lemma}\label{lem:cardinal}
The class of spaces of type at most $\kappa$ is preserved under taking subspaces, closures of subspaces in Hausdorff locally convex spaces, continuous linear images, and finite products. The closed linear span of a family of at most $\kappa$ subspaces of type at most $\kappa$ has type at most $\kappa$.
\end{lemma}
\begin{proof}
We prove the subspace assertion first for normed spaces. Let $X$ have type at most $\kappa$, let $Y\subset X$ be a subspace, and choose $A\subset X$ with $|A|\le\kappa$ and dense linear span. For each finite $H\subset A$ and $n\ge1$ put $F_H=\spn H$, $B_n=\{x:\|x\|<2^{-n}\}$ and
\[
 D_{H,n}=F_H\cap(Y+B_n).
\]
This is an absolutely convex subset of the finite-dimensional space $F_H$, so by \cite[Lemma 3.2.2, p. 90]{PGS} it is the absolutely convex hull of a countable set $G_{H,n}\subset D_{H,n}$. For each $g\in G_{H,n}$ choose $y_g\in Y$ with $\|g-y_g\|<2^{-n}$. Collect all these $y_g$ in a set $A_Y$ of cardinality at most $\kappa$ (choices are indexed by $H,n,g$). For $y\in Y$ and $n\ge1$ choose $v\in\spn A$ with $\|v-y\|<2^{-n}$. Then $v\in D_{H,n}$ for some finite $H$. Expressing $v$ as a finite absolutely convex combination of members of $G_{H,n}$ and replacing each generator by its chosen $y_g$ gives a member of $\spn A_Y$ at distance less than $2^{-n}$ from $y$, by ultrametricity. Thus $\spn A_Y$ is dense in $Y$.

For a locally convex subspace $Y$ of $X$ and a continuous seminorm $q$ on $Y$, there is a continuous seminorm $p$ on $X$ such that $q\le p|_Y$. The image of $Y$ in $X_p$ has type at most $\kappa$ by the normed-space result. Its continuous linear image $Y_q$, under the natural map $y+\ker(p|_Y)\mapsto y+\ker q$, has the same type bound. Continuous images preserve the property since a continuous seminorm on an image pulls back to one on the domain; the corresponding normed quotients are isometric. If a subspace is dense in a larger space, its image in every normed quotient of that larger space is dense, so passage to closure preserves the property. Finite products follow from domination by a maximum of finitely many continuous seminorms. Finally, let $Z$ be the closed linear span of at most $\kappa$ such subspaces. For each continuous seminorm $p$ on $Z$, the images of the subspaces in $Z_p$ have type at most $\kappa$ and together have dense linear span. Choose a generating set of cardinality at most $\kappa$ in each image; the union of at most $\kappa$ such sets has cardinality at most $\kappa$ and dense linear span in $Z_p$.
\end{proof}

For Hausdorff locally convex spaces $E,F$, write
\[
 \mathcal R_\kappa(E,F)=\{S\in\Lop(E,F): S(E)\text{ with its induced topology is of type at most }\kappa\}.
\]
\begin{proposition}\label{prop:cardideal}
The classes $\mathcal R_\kappa(E,F)$ form a linear operator ideal containing the compactoid operators. It is sequentially closed under pointwise convergence within $\Lop(E,F)$; for normed spaces it is in particular closed in operator norm. If $E$ is not of type at most $\kappa$, then $\mathcal R_\kappa(E,E)$ is a proper ideal of $\Lop(E)$.
\end{proposition}
\begin{proof}
For additivity, $(S_1+S_2)(E)$ is a subspace of the image of $S_1(E)\times S_2(E)$ under addition; Lemma~\ref{lem:cardinal} applies. Composition on either side with a continuous linear operator preserves the property by the subspace and image assertions. The ideal contains finite-rank operators, and contains compactoid operators by \cite[Theorem 8.1.3(iv)]{PGS}. If $S_n\to S\in\Lop(E,F)$ pointwise with each $S_n\in\mathcal R_\kappa(E,F)$, then $S(E)\subset\cl{\spn(\bigcup_n S_n(E))}$, a space of type at most $\kappa$ by the last assertion of Lemma~\ref{lem:cardinal}. Properness follows since the identity of $E$ is not in the ideal.
\end{proof}

\begin{theorem}\label{thm:cardinal}
Let $E$ be Hausdorff locally convex and not of type at most $\kappa$, and $T\in\Lop(E)$.
\begin{enumerate}
\item[(i)] If $\{T\}'\cap\mathcal R_\kappa(E,E)$ contains a nonzero operator $S$, then $\cl{S(E)}$ is a nontrivial closed invariant subspace of $T$ of type at most $\kappa$.
\item[(ii)] If $\{T\}''\cap\mathcal R_\kappa(E,E)$ contains a nonzero operator $S$, then this subspace is hyperinvariant for $T$.
\end{enumerate}
In particular, a nonzero operator in $\mathcal R_\kappa(E,E)$ has a nontrivial closed hyperinvariant subspace. No completeness or dual-separation assumption on $E$ is needed.
\end{theorem}
\begin{proof}
By Lemma~\ref{lem:cardinal}, $M=\cl{S(E)}$ has type at most $\kappa$, so $0\ne M\ne E$. Every continuous linear operator commuting with $S$ leaves $M$ invariant. This proves both statements.
\end{proof}

\begin{corollary}\label{cor:cardalg}
Let $E$ be as in Theorem~\ref{thm:cardinal}, and suppose that the continuous dual of $E$ separates points. If the class of $T\in\Lop(E)$ in $\Lop(E)/\mathcal R_\kappa(E,E)$ is algebraic, then $T$ has a nontrivial closed hyperinvariant subspace if and only if no irreducible polynomial in $\K[X]$ annihilates $T$.
\end{corollary}
\begin{proof}
There exists a nonzero polynomial $Q$ with $Q(T)\in\mathcal R_\kappa(E,E)$. If $Q(T)\ne0$, apply Theorem~\ref{thm:cardinal}(ii) to $S=Q(T)$. If $Q(T)=0$, use Proposition~\ref{prop:alg}. In either case the same proposition excludes hyperinvariance when an irreducible polynomial annihilates $T$ (that direction uses the separating-dual assumption).
\end{proof}

\begin{remark}
If $E$ is not of countable type, every $T\in\Lop(E)$ has a nontrivial closed invariant subspace, without any small-range hypothesis: for $x\ne0$ use $\cl{\spn\{T^nx:n\ge0\}}$, which has a countable dense linear span and is thus of countable type (see Section~3). It is proper because $E$ is not of countable type. The preceding criteria are thus of particular interest for hyperinvariance, and the first assertion of Theorem~\ref{thm:cardinal} additionally identifies a subspace from the range of the given commuting operator.
\end{remark}

\begin{corollary}\label{cor:cardc0}
Let $I$ be a set with $|I|>\kappa$ and let $S\in\Lop(c_0(I,\K))$ have nonzero range of type at most $\kappa$. Then the closure of its range $M=\cl{S(c_0(I,\K))}$ is a nontrivial closed hyperinvariant subspace contained in $c_0(J,\K)$ for some $J\subset I$ with $|J|\le\kappa$. Its algebraic codimension is at least $|I|$. For $\kappa=\aleph_0$, for every $\varepsilon>0$ there is a continuous projection onto $M$ of norm at most $1+\varepsilon$.
\end{corollary}
\begin{proof}
Choose at most $\kappa$ vectors whose linear span is dense in $M$. Every vector in $c_0(I,\K)$ has at most countable support, so the union $J$ of their supports has cardinality at most $\kappa$. Hence $M\subset c_0(J,\K)$. This also proves that $c_0(I,\K)$ is not of type at most $\kappa$ if $|I|>\kappa$. The unit vectors indexed by $I\setminus J$ are linearly independent modulo $M$, giving the codimension bound. The hyperinvariance follows from Theorem~\ref{thm:cardinal}. In the countable case, use \cite[Theorem 2.5.9(ii), p. 44]{PGS} for the projection. By Gruson's theorem \cite[Theorem 5.9, pp. 174--176]{vanRooij}, $M$ has an orthogonal basis, which is orthonormal after rescaling because the sup norm takes values in $|\K|$. We do not claim that such a projection commutes with $S$, or claim complementability for general uncountable $\kappa$.
\end{proof}

For example, if $J\subset I$ is infinite and $|J|\le\kappa<|I|$, the coordinate projection onto $c_0(J,\K)$ belongs to $\mathcal R_\kappa(c_0(I,\K),c_0(I,\K))$ but is not compactoid, since its restriction to $c_0(J,\K)$ is the identity. Thus the small-range criterion genuinely goes beyond compactoid operators. The case $\kappa=\aleph_0$ for compactoid operators is an immediate application of the classical theorem on their ranges; for Banach spaces this was noted explicitly in \cite[Remark 10, pp. 212--213]{BE}. It is not a new compactoid spectral result. The ideal need not be closed for arbitrary pointwise convergent nets: the net of coordinate projections onto finite subsets of $I$, directed by inclusion, converges pointwise to the identity on $c_0(I,\K)$, which does not belong to the ideal when $|I|>\kappa$.

\section*{Use of generative AI}
ChatGPT (OpenAI) was used substantially in the exploration of mathematical arguments, proof checking, literature searches, and drafting and revision of this manuscript. The author takes responsibility for the mathematical content, the proofs, and the bibliographic claims of the manuscript.


\begin{thebibliography}{99}
\small\raggedright\setlength{\itemsep}{1pt}\setlength{\parsep}{0pt}
\bibitem{vanRooij} A. C. M. van Rooij, \emph{Non-Archimedean Functional Analysis}, Marcel Dekker, New York, 1978.
\bibitem{PGS} C. P\'erez-Garc\'ia and W. H. Schikhof, \emph{Locally Convex Spaces over Non-Archimedean Valued Fields}, Cambridge Studies in Advanced Mathematics 119, Cambridge University Press, Cambridge, 2010.
\bibitem{Serre} J.-P. Serre, Endomorphismes compl\`etement continus des espaces de Banach $p$-adiques, \emph{Publ. Math. Inst. Hautes \'Etudes Sci.} \textbf{12} (1962), 69--85. \href{https://doi.org/10.1007/BF02684276}{doi:10.1007/BF02684276}.
\bibitem{Gruson} L. Gruson, Th\'eorie de Fredholm $p$-adique, \emph{Bull. Soc. Math. France} \textbf{94} (1966), 67--95. \href{https://www.numdam.org/item/BSMF_1966__94__67_0/}{NUMDAM}.
\bibitem{SliwaFred} W. \'Sliwa, On Fredholm operators between non-Archimedean Fr\'echet spaces, \emph{Compos. Math.} \textbf{139} (2003), 113--118. \href{https://doi.org/10.1023/B:COMP.0000005075.84696.f8}{doi:10.1023/B:COMP.0000005075.84696.f8}.
\bibitem{BE} M. Babahmed and A. El Asri, Invariant subspace problem and compact operators on non-Archimedean Banach spaces, \emph{Extracta Math.} \textbf{35} (2020), 205--219. \href{https://doi.org/10.17398/2605-5686.35.2.205}{doi:10.17398/2605-5686.35.2.205}.
\bibitem{EKB} A. El Asri, A. Kubzdela and M. Babahmed, The Lomonosov type theorems and the invariant subspace problem for non-archimedean Banach spaces, \emph{J. Math. Anal. Appl.} \textbf{543} (2025), 129043. \href{https://doi.org/10.1016/j.jmaa.2024.129043}{doi:10.1016/j.jmaa.2024.129043}.
\bibitem{Benedetto} R. L. Benedetto, \emph{Non-Archimedean Dynamics in Dimension One: Lecture Notes}, Arizona Winter School, 2010. \href{https://swc-math.github.io/aws/2010/2010BenedettoNotes-09Mar.pdf}{Lecture notes}.
\bibitem{Sliwa} W. \'Sliwa, The invariant subspace problem for non-Archimedean Banach spaces, \emph{Canad. Math. Bull.} \textbf{51} (2008), 604--617. \href{https://doi.org/10.4153/CMB-2008-060-9}{doi:10.4153/CMB-2008-060-9}.
\bibitem{Kasprzak} H. Kasprzak, The invariant subspace problem for non-Archimedean K\"othe spaces, \emph{J. Math. Anal. Appl.} \textbf{453} (2017), 1086--1110. \href{https://doi.org/10.1016/j.jmaa.2017.04.034}{doi:10.1016/j.jmaa.2017.04.034}.
\bibitem{KasprzakCorr} H. Kasprzak, Corrigendum to ``The invariant subspace problem for non-Archimedean K\"othe spaces'', \emph{J. Math. Anal. Appl.} \textbf{459} (2018), 1296--1299. \href{https://doi.org/10.1016/j.jmaa.2017.07.071}{doi:10.1016/j.jmaa.2017.07.071}.
\bibitem{Kramar97} E. Kramar, Invariant subspaces for some operators on locally convex spaces, \emph{Comment. Math. Univ. Carolin.} \textbf{38} (1997), 635--644.
\bibitem{Kramar00} E. Kramar, Hyperinvariant subspaces for some operators on locally convex spaces, \emph{Int. J. Math. Math. Sci.} \textbf{23} (2000), 807--814. \href{https://doi.org/10.1155/S0161171200002933}{doi:10.1155/S0161171200002933}.
\bibitem{EB2023} A. El Asri and M. Babahmed, A spectral quasinilpotent operator and the invariant subspace problem for non-Archimedean Banach spaces, \emph{Ann. Funct. Anal.} \textbf{14} (2023), Article 52. \href{https://doi.org/10.1007/s43034-023-00277-x}{doi:10.1007/s43034-023-00277-x}.
\bibitem{Hoover} T. B. Hoover, Hyperinvariant subspaces for $n$-normal operators, \emph{Acta Sci. Math. (Szeged)} \textbf{32} (1971), 109--119. \href{https://acta.bibl.u-szeged.hu/14296/}{Acta archive}.
\end{thebibliography}
\end{document}